\documentclass{article}
\usepackage{amsmath, epsfig,amssymb, multicol}
\newtheorem{lemma}{Lemma}
\newtheorem{theorem}{Theorem}

\newtheorem{corollary}{Corollary}

\newcommand{\G}{{\mathcal G}}
\newcommand{\tr}{{\rm Trace}}
\newcommand{\E}{{\mathrm E}}
\newcommand{\Var}{{\rm Var}}
\newcommand{\Cov}{{\rm Cov}}

\title{Spectra of edge-independent random graphs}
\author{Linyuan Lu
\thanks{University of South Carolina, Columbia, SC 29208,
({\tt lu@math.sc.edu}). This author was supported in part by NSF
grant DMS 1000475. }
  \and Xing Peng
\thanks{University of South Carolina, Columbia, SC 29208,
({\tt pengx@mailbox.sc.edu}).This author was supported  by the Dean's Dissertation Fellowship of College of Arts and Sciences. }}
\begin{document}
\maketitle

\begin{abstract}
  Let $G$ be a random graph on the vertex set $\{1,2,\ldots, n\}$ such
  that edges in $G$ are determined by independent random indicator
  variables, while the probability $p_{ij}$ for $\{i,j\}$ being an
  edge in $G$ is not assumed to be equal. Spectra of the adjacency
  matrix and the normalized Laplacian matrix of $G$ are recently
  studied by Oliveira and Chung-Radcliffe. Let $A$ be the adjacency
  matrix of $G$, $\bar A=\E(A)$, and $\Delta$ be the maximum expected
  degree of $G$.  Oliveira first proved that almost surely $\|A-\bar
  A\|=O(\sqrt{\Delta \ln n})$ provided $\Delta\geq C \ln n$ for some
  constant $C$. Chung-Radcliffe improved the hidden constant in the
  error term using a new Chernoff-type inequality for random matrices.
  Here we prove that almost surely $\|A-\bar A\|\leq
  (2+o(1))\sqrt{\Delta}$ with a slightly stronger condition $\Delta\gg
  \ln^4 n$.  For the Laplacian $L$ of $G$, Oliveira and
  Chung-Radcliffe proved similar results $\|L-\bar L|=O(\sqrt{\ln
    n}/\sqrt{\delta})$ provided the minimum expected degree
  $\delta\gg \ln n$; we also improve their results by removing the
  $\sqrt{\ln n}$ multiplicative factor from the error term under some mild
  conditions. Our results naturally apply to the classic
  Erd\H{o}s-R\'enyi random graphs, random graphs with given expected
  degree sequences, and bond percolation of general graphs.
\end{abstract}

\section{Introduction}
Given an $n \times n$ symmetric matrix $M$, let $\lambda_1(M)$,
$\lambda_2(M),\ldots, \lambda_n(M)$ be the list of eigenvalues of $M$
in the non-decreasing order. What can we say about these eigenvalues
if $M$ is a matrix associated with a random graph $G$? Here $M$ could be
the adjacency matrix (denoted by $A$) or the normalized Laplacian
matrix (denoted by $L$). Both spectra of $A$ and $L$ can be used to
infer  structures of $G$. For example, the spectrum of $A$ is
related to the chromatic number and the independence number.  The
spectrum of $L$ is connected to the mixing-rate of random walks, the
diameters, the neighborhood expansion, the Cheeger constant, the
isoperimetric inequalities, the expander graphs, the quasi-random graphs. For more applications of spectra of the adjacency matrix and the Laplacian matrix, please refer to monographs  \cite{fan4,cv}.


Spectra of adjacency matrices and  normalized Laplacian
matrices of random graphs were extensively investigated in the
literature. For the Erd\H{o}s-R\'enyi random graph model $G(n,p)$,
F\"uredi and Koml{\'o}s \cite{fk} proved that almost surely
$\lambda_n(A)=(1+o(1))np$ and $\max_{i \leq n-1}|\lambda_1(A)|\leq (2+o(1))\sqrt{np(1-p)}$
provided $np(1-p)\gg \log^6 n$; similar results are proved for sparse random graphs \cite{fo,ks} and general random matrices \cite{dj,fk}.
Alon,
Krivelevich, and Vu \cite{akv} showed that with high probability the
$s$-th largest eigenvalue of a random symmetric matrix with
independent random entries of absolute value at most one concentrates
around its median. Chung, Lu, and Vu \cite{flv1,flv2} studied spectra
of adjacency matrices of random power law graphs and spectra of
Laplacian matrices of random graphs with given expected degree
sequences. Their results on random graphs with given expected degree
sequences were complemented by Coja-Oghlan \cite{og,o} for sparser
cases.  For random $d$-regular graphs, Friedman (in a series of papers)
\cite{fks,friedman,friedman2} proved that  the   second largest
eigenvalue (in absolute value) of random $d$-regular graphs is at most  $(2 + o(1))
\sqrt{d-1}$ almost surely  for any $d \geq 4$.

In this paper, we study spectra of the adjacency matrices and the 
 Laplacian matrices of edge-independent random graphs.  Let $G$ be an
edge-independent random graph on the vertex set
$[n]:=\{1,2,\ldots,n\}$; two vertices $i$ and $j$ are adjacent in $G$
with probability $p_{ij}$ independently.  Here $\{p_{ij}\}_{1\leq i,j \leq n}$ are not assumed to be equal.  Let $\bar A:=(p_{ij})_{i,j=1}^n$ be the
expectation of the adjacency matrix $A$ and
$\Delta$ be the maximum expected degree of $G$.  Oliveira \cite{oli}
proved $\|A-\bar A\|=O(\sqrt{\Delta \ln n})$ provided $\Delta\geq C\ln
n$ for some constant $C$.  Chung and Radcliffe \cite{fan7}
improved the hidden constant in the error term  using a Chernoff-type inequality
for random matrices.  We manage to remove the $\sqrt{\ln n}$-factor
from the error term with a slightly stronger assumption on $\Delta$.
We have the following theorem.
\begin{theorem} \label{thm:1}
Consider an edge-independent random graph $G$. 
If $\Delta \gg \ln^4 n$, then almost surely
$$
|\lambda_i(A)-\lambda_i(\bar A) |\leq (2+o(1)) \sqrt{\Delta}
$$
for each $1 \leq i \leq n$.
\end{theorem}

Let $T$ be the diagonal matrix of expected degrees.
Define $\bar L=I - T^{-1/2}\bar A T^{-1/2}$. The matrix $\bar L$ can be viewed
as the ``expected Laplacian'' of $G$.  Oliveira \cite{oli} and Chung-Radcliffe \cite{fan7} proved  theorems   on $L$ which are similar to those on  the adjacency matrix $A$. We are able to improve their results by removing  the $\sqrt{\ln n}$-factor from the error term with some conditions.  
We say that $\bar L$ is  well-approximated by a rank  $k$-matrix if 
there is a $k$ such that all but $k$ eigenvalues
$\lambda_i(\bar L)$ satisfy $|1-\lambda_i(\bar L)|=o(1/\sqrt{\ln n})$.
To make the definition rigorous, 
let  $$\Lambda:=\{\lambda_i(L)\colon  |1-\lambda_i(L)| \geq 1/(g(n)\sqrt{\ln n})\},$$
where $g(n)$ is an arbitrarily slowly growing function; then we have $k:=|\Lambda|$. We have the following theorem.

\begin{theorem} \label{thm:2} Consider an edge-independent random
  graph $G$.  Let $\Lambda$ and $k$ be
defined above. If the minimum expected degree $\delta$ satisfies $\delta\gg  \max\{k, \ln^4 n\}$, then almost surely,
$$
| \lambda_i(L)-\lambda_i(\bar L) |\leq \left(2+\sqrt{\sum_{\lambda\in \Lambda}(1-\lambda)^2}+o(1)\right)\frac{1}{\sqrt{\delta}}$$
for each $1 \leq i \leq n$.
\end{theorem}

Note   ${\rm rank}(I-\bar L)={\rm rank}(\bar A)$. We have the following corollary.
\begin{corollary} \label{cor:1}
 Consider an edge-independent random graph $G$ with ${\rm rank}(\bar A)=k$.
 If the minimum expected degree $\delta$ satisfies $\delta\gg  \max\{k, \ln^4 n\}$, then almost surely,
 we have  $$
| \lambda_i(L)-\lambda_i(\bar L) |\leq \frac{2+\sqrt{k}+o(1)}{\sqrt{\delta}}$$ 
for $1\leq i \leq n$.
\end{corollary}

A special case is the  random graph $G({\mathbf w})$ with given  expected degree sequence 
${\mathbf w}=(w_1, w_2, \ldots,w_n)$, where $v_iv_j$ is an edge with
probability $p_{ij}=\frac{w_iw_j}{\sum_{l=1}^nw_l}.$ 
Let $\delta=w_{\min}$ and 
$\bar w=\frac{\sum_{i=1}^n w_i}{n}$. 
Chung-Lu-Vu  \cite{flv2} proved if  $\delta\gg \sqrt{\bar w} \ln^3 n$ then for each non-trivial eigenvalue $ \lambda_i(L)$, we have
\begin{equation}
  \label{eq:gw}
  |1- \lambda_i(L)|\leq  \frac{4+o(1)}{\sqrt{\delta}}.
\end{equation}

Note that in this case $I-{\bar L}=T^{-1/2}\bar AT^{-1/2}$;
its $(i,j)$-entry is given by $\frac{\sqrt{w_iw_j}}{\sum_{l=1}^nw_l}$.
Thus $I-{\bar L}$ is a rank-1 matrix with non-zero eigenvalues equal 1. 
Hence all non-trivial eigenvalues of $\bar L$ are $1$. Applying  Corollary \ref{cor:1},
we get 
\begin{equation}
  \label{eq:gw2}
  |1- \lambda_i(L)|\leq  \frac{3+o(1)}{\sqrt{\delta}},
\end{equation}
provided $\delta\gg \ln^4 n$.

In comparison to inequality \eqref{eq:gw}, inequality \eqref{eq:gw2} improves the
constant factor with a weaker condition.

Here is another application. 
Let $G$ be a host graph with vertex set $[n]$. The bond percolation of $G$ (with probability $p$) 
is a random spanning subgraph $G_p$ of $G$ such that for each edge $\{i,j\}$ of $G$,
$\{i,j\}$ is retained as an edge of $G_p$ with probability $p$ independently. 
The Erd\H{o}s-R\'enyi graph $G(n,p)$ can be viewed as the  bond percolation
of the complete graph $K_n$. We have the following theorems on the spectrum of $G_p$ for a general graph $G$.

\begin{theorem}
Suppose that the maximum degree $\Delta$ of $G$ satisfies $\Delta\gg \ln^4 n$.
For $p\gg \frac{\ln^4n}{\Delta}$, almost surely we have
$$|\lambda_i(A(G_p))-p\lambda_i(A(G))|\leq (2+o(1))\sqrt{p\Delta}.$$
\end{theorem}

\begin{theorem}
Suppose that all but $k$ Laplacian eigenvalues $\lambda$ of $G$
satisfies $|1-\lambda|=o(\frac{1}{\sqrt{\ln n}})$. 
If the minimum degree $\delta$ of $G$ satisfies $\delta\gg \max\{k,\ln^4 n\}$,
then for $p\gg \max\{\frac{k}{\delta},\frac{\ln^4n}{\delta}\}$, almost surely we have
$$|\lambda_i(L(G_p))-\lambda_i(L(G))|\leq 
\left(2+\sqrt{\sum_{i=1}^k (1-\lambda_i)^2}+o(1))\right)\frac{1}{\sqrt{p\delta}},$$
where $\lambda_1, \ldots, \lambda_k$ are those $k$ Laplacian eigenvalues of $G$ do not satisfy
$|1-\lambda|=o(\frac{1}{\sqrt{\ln n}})$.
\end{theorem}

The rest of the paper is organized as follows. In section 2, we will 
generalize Vu's result  \cite{vu} on the spectral
bound of a random symmetric matrix; we use it to prove theorem \ref{thm:1}. 
In section 3, we will prove several lemmas for  Laplacians. Finally,
 Theorem \ref{thm:2} will be proved in Section 4. 

\section{Spectral bound of random symmetric matrices}
For any matrix $M$, the {\em spectral norm} $\|M\|$ is the largest singular value of $M$; i.e., we have
$$\|M\| =\sqrt{\lambda_{max}(M^*M)}.$$
Here $M^*$ is the conjugate transpose of $M$ and $\lambda_{max}(\bullet)$
is the largest eigenvalue. When $M$ is a symmetric matrix with non-negative entries, we have 
$\|M\|=\lambda_{max}(M).$

We will estimate the spectral norm of 
random symmetric  matrices. Let us start with the following theorem proved by Vu in \cite{vu}.

\begin{theorem} \label{thm:vu}
There are constants $C$ and $C'$ such that
the following holds. Let $b_{  ij}$, $1\leq i\leq j \leq n$ be independent
random variables, each of which has mean $0$ and variance at most $\sigma^2$
and is bounded in absolute value by $K$, where $\sigma\geq C'n^{-1/2}K\ln^2n$.
Then almost surely
\begin{equation}
  \label{eq:vu}
\|B\|\leq 2\sigma\sqrt{n}+C(K\sigma)^{1/2}n^{1/4}\ln n.  
\end{equation}
\end{theorem}

Vu's theorem is already in a very general form; it improves F\"uredi-Koml{\'o}s's result \cite{fk} on $G(n,p)$.  When we consider an edge-independent random graph $G$, let $A$ be the adjacency matrix of $G$ and $\bar A$ be the expectation of $A$. 
If we apply Theorem \ref{thm:vu} to $B:=A-\bar A$, we get
\begin{equation} \label{eq:vu1}
\|A-\bar A\|\leq 2\sigma\sqrt{n}+C(\sigma)^{1/2}n^{1/4}\ln n,
\end{equation}
where $\sigma=\max_{1\leq i\leq j \leq n}\{\sqrt{p_{ij}(1-p_{ij})}\}$.
The upper bound in inequality \eqref{eq:vu1} is weaker than the one in Theorem \ref{thm:1}; this is because the uniform bounds on
$K$ and $\sigma^2$ are too coarse.

To overcome the deficiency, we assume that $b_{ij}$ ($1\leq i \leq j\leq n$) are independent random variables with the
following properties:
\begin{itemize}
\item $|b_{ij}|\leq K$ for $1\leq i<j\leq n$;
\item $\E(b_{ij})=0$, for all $1\leq i<j\leq n$;
\item $\Var(b_{ij})\leq \sigma_{ij}^2$.
\end{itemize}
If $i>j$, we set 
 $b_{ji}=b_{ij}$ and $\sigma_{ji}=\sigma_{ij}$.
Consider a random  symmetric  matrix $B=(b_{ij})_{i,j=1}^n$. 
The following theorem generalizes Vu's theorem.

\begin{theorem}\label{thm:3}
There are constants $C$ and $C'$ such that the following holds.
  Let $B$ be the random symmetric  matrix defined above and 
  $\Delta=\max_{1\leq i\leq n}\sum_{j=1}^n\sigma_{ij}^2$.  If $\Delta\geq C' K^2 \ln^4 n$,
then almost surely 
\begin{equation}
  \label{eq:B}
  \|B\|\leq 2\sqrt{\Delta} + C \sqrt{K} \Delta^{1/4}\ln n.
\end{equation}
\end{theorem}
When $\sigma_{ij}\equiv \sigma$, we have
$\Delta=n\sigma^2$. Thus, inequality \eqref{eq:B} implies
inequality \eqref{eq:vu}. (The condition $\Delta\geq C' K^2 \ln^4 n$
becomes $\sigma\geq \sqrt{C'}n^{-1/2}K\ln^2 n$.)

Replacing $B$ by $cB$, $K$ by $cK$, and $\Delta$ by $c\Delta$,  inequality \eqref{eq:B} is invariant under scaling. Without loss of generality,
we can assume  $K=1$ (by scaling a factor  $\frac{1}{K}$). 
We further assume that  diagonal entries are zeros.
Changing  diagonal entries to zeros can affect the spectral norm by at most $K$,
which is negligible in comparison to the upper bound. 

We use Wigner's trace method \cite{wigner}. We have
\begin{equation}
  \label{eq:tr1}
\E\left(\tr(B^k)\right)=\sum_{i_1,i_2,\ldots, i_k} \E(b_{i_1 i_2}
b_{i_2 i_3} \ldots b_{i_{k-1} i_{k}}b_{i_k i_1}).  
\end{equation}

Each sequence $w:=i_1 i_2 \ldots i_{k-1}i_k i_1$ is a {\it closed
walk} of length $k$ in the complete graph $K_n$. 
Let $E(w)$
 be the set of edges appearing in $w$.
For each edge $e\in E(w)$, let $q_e$ be the number of occurrence
of the edge $e$ in the walk $w$.
By the independence assumption for edges, we can rewrite equation \eqref{eq:tr1} as
\begin{equation}
  \label{eq:tr2}
\E\left(\tr(B^k)\right)=\sum_{w} \prod_{e\in E(w)} \E(b^{q_e}_e).
\end{equation}
Here the summation is taken over all closed walks of length $k$.
If $q_e=1$ for some $e\in E(w)$, then $\prod_{e\in E(w)} \E(c^{q_e}_e)=0$.
Thus we need only to
consider all closed walks such that each edge appears at least twice. 

A closed walk $w$ is {\em good} if each edge in $\E(w)$
occurs more than once. The set of all good closed walks  of length $k$ in $K_n$
is denoted by ${\mathcal G}(n,k)$.

Since $q_e\geq 2$ and $|b_e|\leq 1$, we have
\begin{equation}
  \label{eq:be}
|\E(b_e^{q_e})|\leq  \E(b_e^{2}) =\Var(b_e^2)\leq \sigma_{e}^2.
\end{equation}

Putting equation \eqref{eq:tr2} and inequality \eqref{eq:be} together, we have
\begin{equation}
  \label{eq:bk}
\left |\E\left(\tr(B^k)\right)\right |\leq \sum_{w\in {\G(n,k)}}
 \prod_{e\in E(w)}   \sigma_e^2.
\end{equation}

Let $\G(n,k,p)$ be the set of good closed walks in $K_n$ of length $k$ and 
with $p$ vertices. The key of the trace method is a good estimate  
on $|\G(n,k,p)|$. 
F\"uredi and Koml{\'o}s \cite{fk} proved
\begin{equation}
  \label{eq:fk}
  |\G(n,k,p)|\leq n(n-1)\cdots (n-p+1) \frac{1}{p}{2p-2 \choose p-1}
{k \choose 2p-2} p^{2(k-2p+2)}.
\end{equation}

Let $\tilde \G(k,p)$ be the set of good closed walks $w$ of length $k$ on $[p]$
where vertices first appear in $w$ in the order $1,2,\ldots, p$. It is easy 
to check $|\G(n,k,p)|= n(n-1)\cdots (n-p+1)|\tilde \G(k,p)|$.
The main contribution from Vu's paper \cite{vu} is the following improved bound
(see \cite{vu}, Lemma 4.1):
\begin{equation}
 \label{eq:vub}
|\tilde \G(k,p)| \leq {k \choose 2p-2} 2^{2k-2p+3} p^{k-2p+2}(k-2p+4)^{k-2p+2}.
\end{equation}

We will use this bound to derive the following Lemma.

\begin{lemma} \label{trace}
For any even integer $k$ such that $k^4\leq \frac{\Delta}{32}$, we have
  \begin{equation}
    \label{eq:trace}
    \left |\E\left(\tr(B^k)\right)\right |\leq 2^{k+2}n\Delta^{k/2}.
  \end{equation}
\end{lemma}

\noindent
{\bf Proof:} Let $[n]^{\underline{p}}:=\{(v_1,v_2,\ldots, v_p)\in
[n]^p\colon v_1,v_2,\ldots, v_p \mbox{ are distinct}\}.$ Define
$$\phi\colon {\mathcal G}(n,k,p) \to \tilde {\mathcal G}(k,p) \times
[n]^{\underline{p}}$$ as follows.  For a good closed walk $w=i_1i_2
\ldots i_k i_1 \in {\mathcal G}(n,k,p)$, let $v_1, v_2, \ldots, v_p$
be the list of $p$ vertices in the order as they appear in
$w$. Replacing $v_i$ by $i$ for $1 \leq i \leq p$, we get a good
closed walk $\tilde w\in \tilde {\mathcal G}(k,p)$.  Now we define
$\phi(w)=(\tilde w, (v_1,v_2,\ldots, v_p))$.  Clearly $\phi$ is a
bijection.

For any $w \in \tilde {\mathcal G}(k,p)$, 
we define a rooted tree $T$ (with root $1$)
on the vertex set $[p]$ as follows:
$$
i_ji_{j+1} \in E(T) \ \textrm{if and only if }\ i_{j+1} \not \in \{i_1,
i_2, \ldots,i_j\}.
$$

Equivalently, the edge $i_ji_{j+1} \in E(T)$ if it brings in a new vertex
when it occurs first time. For $2\leq l \leq p$, let $\eta(l)$ be the parent
of $l$. Since $ \sigma^2_{e} \leq K^2=1$ always holds, we can discard
those terms  $ \sigma^2_{e}$ for $e\not\in T(w)$. We get
  \begin{align*}
 \sum_{w \in {\mathcal G}(n,k,p)} \prod_{e\in E(w)} \sigma^2_e&=
 \sum_{\tilde w \in \tilde {\mathcal G}(k,p)} \sum_{(v_1,\ldots, v_p)\in [n]^{\underline{p}}}
 \prod_{xy\in E(\tilde w)} \sigma^2_{v_xv_y} \\
 &\leq  \sum_{\tilde w \in \tilde {\mathcal G}(k,p)}  \sum_{v_1=1}^n \sum_{v_2=1}^n \cdots \sum_{v_{p}=1}^n \prod_{xy\in E(T)} \sigma^2_{v_xv_y}\\
&= \sum_{\tilde w \in \tilde {\mathcal G}(k,p)} \sum_{v_1=1}^n \sum_{v_2=1}^n \cdots \sum_{v_{p}=1}^n
\prod_{y=2}^p \sigma^2_{v_{\eta(y)}v_y}\\
&=  \sum_{\tilde w \in \tilde {\mathcal G}(k,p)} \sum_{v_1=1}^n \sum_{v_2=1}^n \cdots \sum_{v_{p-1}=1}^n
\prod_{y=2}^{p-1} \sigma^2_{v_{\eta(y)}v_y} \sum_{v_p=1}^n  \sigma^2_{v_{\eta(p)}v_p}
\\
&\leq \Delta
\sum_{\tilde w \in \tilde {\mathcal G}(k,p)} \sum_{v_1=1}^n \sum_{v_2=1}^n \cdots \sum_{v_{p-1}=1}^n
\prod_{y=2}^{p-1} \sigma^2_{v_{\eta(y)}v_y}\\
 & \leq \cdots \\
& \leq \Delta^{p-1}  \sum_{\tilde w \in \tilde {\mathcal G}(k,p)}  \sum_{v_1=1}^n 1.\\
&= n \Delta^{p-1}  \left|\tilde {\mathcal G}(k,p)\right|.
\end{align*}
Combining it with inequality \eqref{eq:bk}, we get
\begin{align*}
  \left |\E\left(\tr(B^k)\right)\right | &\leq \sum_{w \in {\mathcal G}(n,k)} \prod_{e\in E(w)} \sigma^2_e \\
&=\sum_{p=2}^{k/2+1}\sum_{w \in {\mathcal G}(n,k)} \prod_{e\in E(w)} \sigma^2_e \\
&\leq \sum_{p=2}^{k/2+1}n \Delta^{p-1} \left|\tilde {\mathcal G}(k,p)\right|\\
&\leq n \sum_{p=2}^{k/2+1} \Delta^{p-1}  {k \choose 2p-2} 2^{2k-2p+3} p^{k-2p+2}(k-2p+4)^{k-2p+2}.
\end{align*}
In the last step, we applied Vu's bound \eqref{eq:vub}.
Let $S(n,k,p):=n \Delta^{p-1}  {k \choose 2p-2} 2^{2k-2p+3} p^{k-2p+2}(k-2p+4)^{k-2p+2}$.
One can show
$$S(n,k,p-1)\leq \frac{16 k^4}{\Delta}S(n,k,p).$$
When $k^4\leq \frac{\Delta}{32}$, we have $S(n,k,p-1)\leq \frac{1}{2}S(n,k,p)$.
Thus,
\begin{align*}
 \left |\E\left(\tr(B^k)\right)\right |
&\leq  \sum_{p=2}^{k/2+1}S(n,k,p)\\
&\leq S(n,k,k/2+1) \sum_{p=2}^{k/2+1} \left(\frac{1}{2}\right)^{k/2+1-p}\\
&<  2S(n,k,k/2+1)\\
&=n 2^{k+2}\Delta^{k/2}.
\end{align*}
The proof of this Lemma is finished. \hfill $\square$

Now we are ready to prove Theorem \ref{thm:3}.

\noindent
{\bf Proof of Theorem \ref{thm:3}:} 
We have
\begin{align*}
\Pr(\|B\|\geq  2\sqrt{\Delta} + C  \Delta^{1/4}\ln n)
&= \Pr(\|B\|^k\geq  (2\sqrt{\Delta} + C  \Delta^{1/4}\ln n)^k)   \\
&\leq \Pr(\tr(\|B\|^k) \geq  (2\sqrt{\Delta} + C \Delta^{1/4}\ln n)^k)   \\
&\leq \frac{\E(\tr(\|B\|^k))}{ (2\sqrt{\Delta} + C  \Delta^{1/4}\ln n))^k}\\
&\leq  \frac{n2^{k+2}\Delta^{k/2}}{ (2\sqrt{\Delta} + C \Delta^{1/4}\ln n))^k}\\
&= 4n e^{-(1+o(1))\frac{C}{2}k\Delta^{-1/4}\ln n}.
\end{align*}
Setting $k= \left(\frac{\Delta}{32}\right)^{1/4}$, we get
$$\Pr(\|B\|\geq  2\sqrt{\Delta} + C  \Delta^{1/4}\ln n)=o(1)$$
for sufficiently large $C$.
The proof of Theorem \ref{thm:3} is finished. \hfill $\square$

\noindent
{\bf Proof of Theorem \ref{thm:1}:} 
Let $B=A-\bar A$. Notice $|b_{ij}|\leq 1$ and
$\Var(b_{ij})=p_{ij}(1-p_{ij})\leq p_{ij}.$

Apply Theorem \ref{thm:3} to $B$ with $K=1$, $\sigma_{ij}^2=p_{ij}$,
and $\Delta=\max_{1\leq i\leq 1} \sum_{j=1}^n p_{ij}$. We get
$$\|B\|\leq 2\sqrt{\Delta}+C  \Delta^{1/4}\ln n.$$
When $\Delta\gg \ln^4 n$, we have
$$\|B\|\leq (2+o(1))\sqrt{\Delta}.$$
Applying Weyl's Theorem, we get
$$|\lambda_i(A)-\lambda_i(\bar A)|\leq \|A-\bar A\|\leq  (2+o(1))\sqrt{\Delta}.$$
The proof of Theorem \ref{thm:1} is completed. \hfill $\square$

\section{Lemmas for Laplacian eigenvalues}


In this section, we will present some necessary lemmas  for proving Theorem \ref{thm:2}. 
Recall $G$ is an edge-independent random graph over
$[n]$ such that  $\{i,j\}$ forms as an edge
with probability $p_{ij}$ independently.  For $1\leq i\leq j \leq n$, let
$X_{ij}$ be the random indicator variable for $\{i,j\}$ being an edge;
we have $\Pr(X_{ij}=1)=p_{ij}$ and $\Pr(X_{ij}=0)=1-p_{ij}$.  For each
vertex $i \in [n]$, we use $d_i$ and $t_i$ to denote the degree and
the expected degree of vertex $i$ in $G$ respectively. 
We have
$$d_i=\sum_{j=1}^n X_{ij} \ \mbox{and} \ t_i=\sum_{j=1}^n p_{ij}.
$$
  Let $D$ (and $T$) be the
diagonal matrix with $D_{ii}=d_i$  (and $T_{ii}=t_i$) respectively.
The matrix $T$ is the expectation of $D$.
Note that  we use $A$ and $L$ to denote
the adjacency matrix and the Laplacian matrix of $G$. Here
$L=I-D^{-1/2}AD^{-1/2}$.  Moreover, we let $\bar A:=\E(A)=(p_{ij})_{i,j=1}^n$ be the
expectation of $A$.
We also define $\bar L=I-T^{-1/2}\bar AT^{-1/2}$. The matrix $\bar L$ can be viewed
as the ``expected Laplacian matrix'' of $G$.

For notational convenience,   we write  eigenvalues of the expected Laplacian matrix $\bar L$ as  
$\mu_1, \mu_2, \ldots, \mu_n$ 
such that
$$|1-\mu_1| \geq |1-\mu_2|\geq \cdots \geq |1-\mu_n|.$$
By the definition of $k$ and $\Lambda$, we have
$\Lambda=\{\mu_1, \ldots, \mu_k\}$ and $|1-
\mu_i|=o(1/\sqrt{\ln n})$ for all $i \geq k+1$. 

For $1 \leq i\leq n$, let $\phi^{(i)}$ be the orthonormal eigenvector  of $\mu_i$ for  $\bar L$. Observe that $\phi^{(i)}$ can also be viewed as the orthonormal
 eigenvector for  $T^{-1/2}\bar AT^{-1/2}$ corresponding to the
 eigenvalue $1-\mu_i$.
 Thus we can rewrite $T^{-1/2}\bar AT^{-1/2}=\sum_{i=1}^n (1-\mu_i)\phi^{(i)} \phi^{(i)'}$.
Let $M= \sum_{i=1}^k (1-\mu_i)\phi^{(i)} \phi^{(i)'}$ and $N=
\sum_{i=k+1}^n (1-\mu_i)\phi^{(i)} \phi^{(i)'}$.  
Observe a fact  $\|N\|=o(1/\sqrt{\ln n})$ and this fact will be used later.  For
a square matrix $B$, we define 
$$f(B)=D^{-1/2}T^{1/2}BT^{1/2}D^{-1/2}-B.$$
We shall rewrite $L-\bar L$ as a sum of four matrices. Notice  $L-\bar L=D^{-1/2} A D^{-1/2}- T^{-1/2} \bar A T^{-1/2}$.
 It is easy to verify $L-\bar L=M_1+M_2+M_3+M_4$, where $M_i$  are following.
\begin{align*}
M_1&=T^{-1/2}(A-\bar A)T^{-1/2},\\
M_2 &= f(M_1),\\
M_3&= f(N),\\
M_4&=f(M).
\end{align*}
Here the matrices $M$ and $N$ are defined above. We
will bound $\| M_i\|$ for $1\leq i \leq 4$ separately.

\begin{lemma}\label{l:M1}
If $\delta\gg \ln^4 n$, then 
$$\|M_1\| \leq (2+o(1))\frac{1}{\sqrt{\delta}}.$$
\end{lemma}
{\bf Proof:} We are going to apply Theorem  \ref{thm:3} to $B=M_1=T^{-1/2}(A-\bar A)T^{-1/2}$.
Note  $$|b_{ij}|\leq (t_it_j)^{-1/2}\leq 1/\delta,$$
and 
$$\Var(b_{ij})=\Var( (t_it_j)^{-1/2} (a_{ij}-p_{ij}))=\frac{p_{ij}(1-p_{ij})}{t_it_j}
\leq \frac{p_{ij}}{t_it_j}.$$
Let $K=1/\delta$ and $\sigma_{ij}^2=  \frac{p_{ij}}{t_it_j}$.
We have 
$$\Delta(B)=\max_{1\leq i \leq n}\sum_{j=1}^n   \frac{p_{ij}}{t_it_j}
\leq \max_{1\leq i \leq n} \sum_{j=1}^n  \frac{p_{ij}}{t_i\delta}
=\frac{1}{\delta}.$$
By  Theorem  \ref{thm:3}, we get
$$\|B\|\leq \frac{2}{\sqrt{\delta}} + C \delta^{-3/4} \ln n.$$
When $\delta\gg \ln^4 n$,  we have $C\delta^{-3/4} \ln n=o(\frac{2}{\sqrt{\delta}})$.
Thus, we get $\|M_1\|\leq \frac{2+o(1)}{\sqrt{\delta}}$ and the proof of the lemma is finished. \hfill $\square$

We have the following lemma on the function $f$.

\begin{lemma}\label{l:lm8}
If $\|B\|=o(1/\sqrt{\ln n})$, then $\| f(B)\|=o(1/\sqrt{\delta})$.
\end{lemma}

Before we prove Lemma \ref{l:lm8}, we have two corollaries.

\begin{corollary} \label{l:M2}
If $\delta\gg \ln^4 n$, then
we have $\| M_2 \|=o(1/\sqrt{\delta}).$
\end{corollary}

\noindent
{\bf Proof:} Recall $M_2=f(M_1)$. By Lemma \ref{l:M1}, we have 
$$\|M_1\|\leq \frac{2+o(1)}{\sqrt{\delta}}=o(1/\sqrt{\ln n}).$$
By Lemma \ref{l:lm8}, we have this  Corollary. \hfill $\square$

\begin{corollary} \label{l:M3}
If $\delta\gg \ln^4 n$, then
we have  $\| M_3 \|=o(1/\sqrt{\delta}).$
\end{corollary}
{\bf Proof:} Recall $M_3=f(N)$. By the definition of $N$, we have $\|N\|=o(1/\sqrt{\ln n})$. Lemma \ref{l:lm8} gives us the Corollary. \hfill $\square$

To prove Lemma \ref{l:lm8}, we need the following Chernoff inequality.
\begin{theorem}{\cite{chernoff}}
\label{chernoff}
 Let $X_1,\ldots,X_n$ be independent random variables with
$$\Pr(X_i=1)=p_i, \qquad \Pr(X_i=0)=1-p_i.$$
We consider the sum $X=\sum_{i=1}^n X_i$,
with expectation $\E(X)=\sum_{i=1}^n p_i$. Then we have
\begin{align*}
\mbox{(Lower tail)~~~~~~~~~~~~~~~~~}
\qquad \qquad  \Pr(X \leq \E(X)-\lambda)&\leq e^{-\lambda^2/2\E(X)},\\
\mbox{(Upper tail)~~~~~~~~~~~~~~~~~}
\qquad \qquad
\Pr(X \geq \E(X)+\lambda)&\leq e^{-\frac{\lambda^2}{2(\E(X) + \lambda/3)}}.
\end{align*}
\end{theorem}

We can use the lemma above to prove the degree of each vertex concentrates around its expectation.

\begin{lemma}\label{l:degree}
Assume $t_i \geq \ln n$ for $1\leq i \leq n$. Then with probability
at least $1-\frac{1}{n^2}$, for all $1 \leq i \leq n$ we have
$$
|d_i-t_i| \leq  3 \sqrt{t_i \ln n}.
$$
\end{lemma}
{\bf Proof:} Recall $X_{ij}$ is the random indicator variable for $\{i,j\}$ being an edge. Note $d_i=\sum_{j=1}^n X_{ij}$ and $\E(d_i)=\sum_{j=1}^n p_{ij}=t_i$.
Applying the lower tail of Chernoff's inequality with $\lambda=3\sqrt{t_i\log n}$,
we have
$$\Pr \left (d_i-t_i\leq - \lambda \right ) \leq e^{-\lambda^2/2t_i}=\frac{1}{n^{9/2}}.$$

Applying the upper tail of Chernoff's inequality with $\lambda=3\sqrt{t_i\log n}$,
we have
$$\Pr \left(d_i-t_i \geq  \lambda \right )
\leq e^{-\frac{\lambda^2}{2(t_i +
 \lambda/3)}}
\leq\frac{1}{n^{27/8}}.$$
The union bound gives the lemma.
\hfill $\square$

By  Lemma \ref{l:degree}, we can write $d_i=(1+o(1))t_i$  for $1 \leq i \leq n$.

\begin{lemma} \label{l:lm7}
When $\delta \gg \ln n$, we have
$$
  \|D^{-1/2}T^{1/2}-I \| =O\left(\sqrt{\frac{\ln n}{ \delta}}
\right) \hspace{0.2cm}  \mbox{and} \hspace{0.2cm} 
\|T^{1/2}D^{-1/2}\|=1+o(1).$$
\end{lemma}
{\bf Proof:} We note that $D^{-1/2}T^{1/2}-I $ is diagonal and
the $(i,i)$-the entry is
$\sqrt{t_i/d_i}-1$.  We have
\begin{align*}
\left|
\frac{\sqrt{t_i}}{\sqrt{d_i}}-1\right|&=\left|\frac{t_i-d_i}{\sqrt{d_i}(\sqrt{t_i}+\sqrt{d_i})}\right| \\
                                      & \leq \left(\frac{3}{2}+o(1)\right) \sqrt{\frac{\ln n}{t_i}}   \\
                                      & = O( \sqrt{\frac{\ln n}{\delta}}).
 \end{align*}
The first part of this lemma is proved while the second part follows from the triangle inequality.
The proof of the lemma is finished. \hfill $\square$

We are ready to prove  Lemma \ref{l:lm8}.
\hspace{0.2cm}

\noindent
{\bf Proof of Lemma \ref{l:lm8}:} Recall that $f(B)=D^{-1/2}T^{1/2}BT^{1/2}D^{-1/2}-B$.  We have
\begin{align*}
f(B)&=D^{-1/2}T^{1/2}BT^{1/2}D^{-1/2}-B\\
    &=D^{-1/2}T^{1/2}BT^{1/2}D^{-1/2}-BT^{1/2}D^{-1/2}+BT^{1/2}D^{-1/2}-B\\
    &=(D^{-1/2}T^{1/2}-I)BT^{1/2}D^{-1/2}+B(T^{1/2}D^{-1/2}-I).
\end{align*}
Recall Lemma \ref{l:lm7}. By the triangle inequality, we have
\begin{align*}
\|f(B)\| &\leq  \|D^{-1/2}T^{1/2}-I\| \|B\| \|T^{1/2}D^{-1/2}\|+\|B\|\|(T^{1/2}D^{-1/2}-I)\|\\
         & \leq O\left(\sqrt{\frac{\ln n}{\delta}}\right)\|B\|(1+o(1))+\|B\| O\left(\sqrt{\frac{\ln n}{\delta}}\right)\\
         &=o\left( \frac{1}{\sqrt{\delta}}\right).
\end{align*}
We use the assumption $\|B\|=o(1/\sqrt{\ln n})$ in the last step and  we completed the proof of the lemma.
\hfill $\square$

\section{Proof of Theorem \ref{thm:2}}
It remains to estimate $\|M_4\|$. Recall $M_4=f(M)$ and 
$M=\sum_{i=1}^k (1-\mu_i) \phi^{(i)} \phi^{(i)'}$.

For $1 \leq i \leq n$, write $\phi^{(i)}$ as a vector 
$(\phi^{(i)}_1,\phi^{(i)}_2,\cdots,\phi^{(i)}_n)'$.  Let $\|\phi^{(i)}\|_\infty$
be the maximum over $\{|\phi^{(i)}_1|,|\phi^{(i)}_2|,\cdots,|\phi^{(i)}_n|\}$.
We have the following lemma.
\begin{lemma} \label{l:eigen}
For each $1 \leq i \leq n$, we have
$$
|1-\mu_i| \cdot \|\phi^{(i)}\|_\infty \leq \frac{1}{\sqrt{\delta}}.
$$
\end{lemma}
{\bf Proof:} Assume $ \|\phi^{(i)}\|_\infty=|\phi^{(i)}_j|$ for some index $j$. Since
$\phi^{(i)}$ is the orthonormal eigenvector associated with $1-\mu_i$
for $T^{-1/2}\bar AT^{-1/2}$, we have $T^{-1/2}\bar AT^{-1/2}
\phi^{(i)}=(1-\mu_i) \phi^{(i)}$. In particular, $(T^{-1/2} \bar AT^{-1/2}
\phi^{(i)})_j=(1-\mu_i) \phi^{(i)}_j$ holds. We have
\begin{align*}
|1-\mu_i| \cdot \|\phi^{(i)}\|_\infty  &= |1-\mu_i| |\phi^{(i)}_j|\\
&=|(T^{-1/2} \bar AT^{-1/2} \phi^{(i)})_j|\\
& \leq \sum_{l=1}^n \frac{p_{jl} |\phi^{(i)}_l| }{\sqrt{t_jt_l}} \\
                             &\leq  \left(\sum_{l=1}^n (\phi^{(i)}_l)^2\right)^{1/2} \left(\sum_{l=1}^n \frac{p_{jl}^2}{t_jt_l}\right)^{1/2}\\
                             & \leq \frac{1}{\sqrt{\delta}}
 \left(\sum_{l=1}^n \frac{p_{jl}}{t_j}\right)^{1/2}\\
                             &= \frac{1}{\sqrt{\delta}}.
\end{align*}
The lemma is proved. \hfill $\square$
\begin{lemma} \label{l:variance}
Assume $\delta \gg \ln n$.
For $1 \leq i \leq n$, consider a random variable $X_i:=\frac{(d_i-t_i)^2}{t_i}$. 
We have $\E(X_i)\leq 1$ and 
$\Var(X_i) \leq 2+o(1)$ for $1 \leq i \leq n$, and  $\Cov(X_i,X_j)=\frac{p_{ij}(1-p_{ij})(1-2p_{ij})}{t_it_j}$ for $1 \leq i \not = j \leq n$.
\end{lemma}
{\bf Proof:} For $1 \leq i< j \leq n$, recall that $X_{ij}$ is the
random indicator  variable for $\{i,j\}$ being an edge.  We define
$Y_{ij}=X_{ij}-p_{ij}$.  Thus we have $d_i-t_i=\sum_{j=1}^n Y_{ij}$.
Note that $\E(Y_{ij})=0$ and $\Var(Y_{ij})=p_{ij}(1-p_{ij})$. We get
$\E(X_i)=\frac{1}{t_i}\sum_{j=1}^{n} p_{ij}(1-p_{ij})\leq 1$.  We have
$$\E(X_i^2)=\frac{1}{t_i^2} \E\left(\sum_{j_1,j_2,j_3,j_4} Y_{ij_1} Y_{ij_2} Y_{ij_3}Y_{ij_4}\right).$$
Since we have $\E(Y_{ij})=0$, the non-zero term occurs either
$j_1=j_2=j_3=j_4$, or $j_1=j_2 \not = j_3=j_4$, or $j_1=j_3\not
=j_2=j_4$, or $j_1=j_4 \not =j_2=j_3$.  The contribution from the
first case is
$$\frac{1}{t_i^2} \sum_{j=1}^n \E(Y_{ij}^4)=\frac{1}{t_i^2} \sum_{j=1}^n (1-p_{ij})^4 p_{ij}+p_{ij}^4(1-p_{ij}) \leq \frac{1}{t_i^2} \sum_{j=1}^n p_{ij}=\frac{1}{t_i}=o(1)$$
as we assume $\delta \gg \ln n$.  The contribution from the second
case is
$$
\frac{1}{t_i^2} \sum_{j_1 \not = j_2} \E(Y_{ij_1}^2) \E(Y_{ij_2}^2)=\frac{1}{t_i^2} \sum_{j_1\not =j_2} p_{ij_1}p_{ij_2}(1-p_{ij_1})(1-p_{ij_2}).
$$
The contribution from the third case and the forth case equal  the contribution from the second case. Thus
$$
\E(X_i^2)=o(1)+\frac{3}{t_i^2} \sum_{j_1 \not =j_2} p_{ij_1}p_{ij_2}(1-p_{ij_1})(1-p_{ij_2}).
$$
We have
\begin{align*}
\E(X_i)^2&=\frac{1}{t_i^2}\sum_{j=1}^n p_{ij}^2(1-p_{ij})^2+\frac{1}{t_i^2}\sum_{j_1 \not =j_2} p_{ij_1}p_{ij_2}(1-p_{ij_1})(1-p_{ij_2})\\
          &=o(1)+ \frac{1}{t_i^2}\sum_{j_1 \not =j_2} p_{ij_1}p_{ij_2}(1-p_{ij_1})(1-p_{ij_2}).
\end{align*}
Therefore we get
\begin{align*}
\Var(X_i)&=\E(X_i^2)-\E(X_i)^2\\
         &=\frac{2}{t_i^2} \sum_{j_1 \not =j_2} \left(p_{ij_1}p_{ij_2}(1-p_{ij_1})(1-p_{ij_2})\right)+o(1) \\
         &\leq  \frac{2}{t_i^2} \sum_{j_1 \not =j_2} \left(p_{ij_1}p_{ij_2}\right)+o(1)\\
         & \leq 2+o(1).
\end{align*}
The covariance $\Cov(X_i,X_j)$ can be computed similarly. Here we omit the details.
 \hfill $\square$

\begin{lemma} \label{l:comb}
For any non-negative numbers $a_1, a_2,\ldots a_n$,  
let $ X=\sum_{i=1}^n a_i \frac{(d_i-t_i)^2}{t_i}$.  Set $a=\max_{1\leq i\leq n}\{a_i\}$.
Then  we have 
$$\Pr\left(X \geq \sum_{i=1}^n a_i + \eta \sqrt{a\sum_{i=1}^n a_i} \ \right)\leq \frac{2+o(1)}
{\eta^2}.$$
\end{lemma}
{\bf Proof:} We have $X=\sum_{i=1}^n a_i X_i$, where
$X_i=\frac{(d_i-t_i)^2}{t_i}$.  We shall use the second moment method
to show that $X$ concentrates around its expectation. By Lemma \ref{l:variance}, 
we have $\E(X_i)\leq 1$. Thus $\E(X)\leq \sum_{i=1}^n a_i$.
For the  variance, by Lemma \ref{l:variance}, we obtain
\begin{align*}
\Var(X) &=  \sum_{i=1}^n a_i^2\Var(X_i)+\sum_{i \not= j}
a_ia_j\Cov(X_i,X_j)\\
& \leq  \sum_{i=1}^n a_i^2 (2+o(1))+  \sum_{i=1}^n \sum_{j=1}^n \frac{a_ia_j p_{ij}(1-p_{ij})(1-2p_{ij})}{t_it_j}\\
        & \leq    (2+o(1)) a \sum_{i=1}^n a_i+ \frac{a}{\delta}\sum_{i=1}^n a_i\sum_{j=1}^n \frac{  p_{ij}}{t_i}\\
 &= (2+o(1))a \sum_{i=1}^n a_i.
\end{align*}
Applying Chebychev's inequality, we have
\begin{align*}
 \Pr\left(X \geq \sum_{i=1}^n a_i + \eta \sqrt{a\sum_{i=1}^n a_i} \hspace{0.2cm}\right)
&\leq \Pr\left(X-\E(X)\geq \eta \sqrt{a\sum_{i=1}^n a_i}\hspace{0.2cm} \right)\\
&\leq \frac{\Var(X)}{\eta^2 a\sum_{i=1}^n a_i}\\
&\leq \frac{2+o(1)}
{\eta^2}.
\end{align*}
\hfill $\square$

We are ready to prove an upper bound on $\|M_4\|$.
\begin{lemma}\label{l:M4} If $\delta\gg \{k, \ln^4 n\}$, then 
we have
$$\|M_4 \|\leq  (1+o(1))) \frac{\sqrt{\mathstrut \sum_{\lambda\in \Lambda} \mathstrut (1-\lambda)^2}}{\sqrt{\delta}}.$$

\end{lemma}
{\bf Proof:} 
Let $\Phi:=(\phi^{(1)}, \ldots,\phi^{(k)})$ be an $n\times k$  matrix such that its columns are mutually orthogonal
and $Q$ be a diagonal $k\times k$ matrix such that $Q_{ii}=1-\mu_i$.
We have $M=\sum_{i=1}^k (1-\mu_i) \phi^{(i)} \phi^{(i)'}=\Phi Q\Phi'$. Thus,
\begin{align*}
M_4& =D^{-1/2}T^{1/2} M T^{1/2} D^{-1/2}-M\\
    &= D^{-1/2}T^{1/2} M T^{1/2} D^{-1/2}-M T^{1/2} D^{-1/2}+M T^{1/2} D^{-1/2} -M \\
    &=(D^{-1/2}T^{1/2} -I) M T^{1/2} D^{-1/2} + M ( T^{1/2}D^{-1/2} - I)\\
&= (D^{-1/2}T^{1/2} -I)\Phi Q \Phi' T^{1/2} D^{-1/2} + \Phi Q \Phi' (T^{1/2} D^{-1/2} - I).
\end{align*}

Let $U=(D^{-1/2}T^{1/2} -I)\Phi Q$.  By the definition of  $\Phi$,
we have $\|\Phi\|=1$. Since $\|T^{1/2} D^{-1/2}\|=1+o(1)$, by the triangle inequality, 
we get
\begin{align*}
  \|M_4\| &= \|U \Phi' T^{1/2} D^{-1/2} + \Phi U'\|\\
&\leq  \|U\| \|\Phi'\| \|T^{1/2} D^{-1/2}\| + \|\Phi\| \| U'\|\\
&=(2+o(1)) \|U\|. 
\end{align*}

By the definition of the norm of a non-square matrix, we have
\begin{align*}
  \|U\| &= \sqrt{\|U U'\|} \\
&\leq \sqrt{\tr(UU')}\\
&= \sqrt{\tr\left((D^{-1/2}T^{1/2}-I)\Phi Q Q \Phi' (T^{1/2}D^{-1/2}-I)\right)}\\
&= \sqrt{\sum_{j=1}^n \sum_{i=1}^k (1-\mu_i)^2\left(\phi^{(i)}_j\right )^2 
\left( \frac{\sqrt{t_j}}{\sqrt{d_j}} -1 \right)^2}.
\end{align*}

Let $a_j:=\sum_{i=1}^k (1-\mu_i)^2\left(\phi^{(i)}_j\right )^2 $.
We have the following estimate on the norm of $U$,
\begin{align*}
\|U\|^2&\leq \sum_{j=1}^n a_j \left( \frac{\sqrt{t_j}}{\sqrt{d_j}} -1 \right)^2\\
&=  \sum_{j=1}^n a_j \frac{(t_j-d_j)^2}{d_j(\sqrt{t_j}+ \sqrt{d_j})^2}\\
&= (1+o(1))  \sum_{j=1}^n a_j \frac{(t_j-d_j)^2}{4t_j^2}\\
&\leq \frac{1+o(1)}{4\delta}  \sum_{j=1}^n a_j \frac{(t_j-d_j)^2}{t_j}.
\end{align*}

 Note that
\begin{align*}
\sum_{j=1}^n a_j &= \sum_{j=1}^n
\sum_{i=1}^k (1-\mu_i)^2\left(\phi^{(i)}_j\right )^2  \\
&=  \sum_{i=1}^k (1-\mu_i)^2\sum_{j=1}^n\left(\phi^{(i)}_j\right )^2  \\
&=\sum_{i=1}^k (1-\mu_i)^2.
\end{align*}
By Lemma \ref{l:eigen}, we have
$$|1-\mu_i| \cdot \|\phi^{(i)}\|_\infty \leq \frac{1}{\sqrt{\delta}}.$$
Hence, for $1\leq j\leq n$, we get 
$$a_j=
\sum_{i=1}^k (1-\mu_i)^2\left(\phi^{(i)}_j\right )^2  \leq 
\frac{k}{\delta}.$$
If we let $a=\max_{1\leq j\leq n}\{a_j\}$, then we have $ a \leq  \frac{k}{\delta}$.

Choose $\eta:=\sqrt[3]{\delta}/\sqrt[3]{k}$; we have $\eta\to \infty$ as 
$n$ approaches the infinity. 
Applying Lemma \ref{l:comb}, 
with probability $1-o(1)$, we have
\begin{align*}
\sum_{j=1}^n a_j \frac{(t_j-d_j)^2}{t_j}
&\leq \sum_{j=1}^n a_j + \eta \sqrt{a \sum_{j=1}^n a_j} \\
&\leq \sum_{i=1}^k (1-\mu_i)^2 + \frac{\sqrt[3]{\delta}}{\sqrt[3]{k}} 
\sqrt{\frac{k}{\delta} \sum_{i=1}^k (1-\mu_i)^2}\\
&= \left(\sqrt{\sum_{i=1}^k (1-\mu_i)^2}+o(1)\right)^2.
\end{align*}
Therefore, we get the following upper bounds on $\|U\|$ and $\|M_4|$;  
$$\|U\|\leq (1+o(1)) \frac{\sqrt{\sum_{i=1}^k (1-\mu_i)^2}}{2\sqrt{\delta}},$$
and 
$$\|M_4\| \leq (1+o(1)) \frac{\sqrt{\sum_{i=1}^k (1-\mu_i)^2}}{\sqrt{\delta}}=(1+o(1))) \frac{\sqrt{\mathstrut \sum_{\lambda\in \Lambda} \mathstrut (1-\lambda)^2}}{\sqrt{\delta}}.$$
We proved the lemma. \hfill $\square$

\noindent {\bf Proof of Theorem \ref{thm:2}:} 
Recall $L-\bar L=M_1+M_2+M_3+M_4$. By the triangle inequality, we have $\|L-\bar L\| \leq
\|M_1\|+\|M_2\|+\|M_3\|+\|M_4\|$. 
Combining  Lemma \ref{l:M1}, Corollary \ref{l:M2}, Corollary \ref{l:M3}, 
and Lemma \ref{l:M4}, we get
\begin{align*}
\|L-\bar L\|&\leq \|M_1\|+\|M_2\|+\|M_3\|+\|M_4\| \\
  &\leq  (2+o(1))\frac{1}{\sqrt{\delta}}+ o\left( \frac{1}{\sqrt{\delta}}\right)+o\left( \frac{1}{\sqrt{\delta}}\right)+  (1+o(1)) \frac{\sqrt{\sum_{i=1}^k (1-\mu_i)^2}}{\sqrt{\delta}}\\
 &=  \left(2+\sqrt{\sum_{i=1}^k (1-\mu_i)^2} +o(1)\right)\frac{1}{\sqrt{\delta}}.
\end{align*}
Finally we apply  Weyl's Theorem. The proof of Theorem \ref{thm:2} is finished. 
 \hfill $\square$

\end{document}